\documentclass[11pt,reqno]{amsart}
\usepackage{graphicx}

\begin{document}
\newtheorem{thm}{Theorem}
 \newtheorem{lem}{Lemma}
 \newtheorem{prop}{Proposition}[section]
 \newtheorem{defn}{Definition}[section]
  \newtheorem{rem}{Remark}
\newtheorem{cor}{Corollary}


%
\title[Some analytically solvable problems of mean-field games]{Some analytically solvable problems of the mean-field games theory}


\author{Sergey I. Nikulin, Olga S. Rozanova}

\address{
              Moscow State University, Moscow, 119991, Russia}
\email
{rozanova@mech.math.msu.su, sergey.nikuline@mail.ru}           







\begin{abstract}
We study the mean field games equations  \cite{1}, consisting of the coupled Kolmogorov-Fokker-Planck and Hamilton-Jacobi-Bellman equations. The equations are complemented by initial and terminal conditions. It is shown that with some specific choice of data, this problem can be reduced to solving a quadratically nonlinear system of ODEs. This situation occurs naturally in economic applications. As an example, the problem of forming an investorТs opinion on an asset is considered.
\end{abstract}

\maketitle


\section{ Introduction}

The theory of Mean field games (MFG) studies the models with a large number of small components (agents) that interact with each other achieving their individual objectives. 
The term "mean field" means that the strategy of each agent to achieve the maximum of its individual utility directly depends on the average distribution of influences of other agents and does not depend on the initial configuration of the system.
 The mean field theory is well known in statistical physics, but similar concepts related to active objects  \cite{1}, \cite{2} were formulated only in the last decade.  These concepts, together with the optimal control theory, have made it possible to study models in economics and sociology.  Now MFG are widely used  in areas requiring analysis of differential games with a large number of participants.

Here is a heuristic derivation of MFG equations, following  \cite{1}. More details can be found in \cite{3}. Let us assume that averaged behavior of the agents is described by Ito stochastic process  $X_s\in \mathbb R$ that is given as
\begin{equation}\label{X}
dX_s = \alpha_s\, d s + \delta\, d W_s, \quad X_t=x,\end{equation}
where $x\in \mathbb R$  is a point in the space of the states, $0\le t\le
s\le T$, $W_s$ is a standard  Brownian motion, $\delta$
is a positive constant,  $\alpha_s$ is a
parameter, choosing the value of which from a given Borel set $U\subset \mathbb R $ at any time $ s $, one can control the process
$X_s$ (the stochastic process $\alpha_s=\alpha(s,X_s)$ is an admissible Markov control). Fixing
$x$, we mark trajectory for a specific agent.

 The problem of control is to define maximum over all the admissible controls $\alpha_t$ of the expresssion
$$J(t,x;\alpha_t) = {\mathbb E}\Bigg[
\int_t^T F(s,X_s;\alpha_s)ds + K(X_T) \Bigg], $$  where $F:{\mathbb
R}_+\times {\mathbb R}\times U\to \mathbb R$ and $K:{\mathbb R}\to
\mathbb R$  are prescribed continuous functions and the process  $X_s$ obeys \eqref{X}.

Let us consider the payoff  function $\Phi: {\mathbb R}_+\times
\mathbb R$ as
$$\Phi (t,x)=\underset{\alpha\in U}{\rm sup}\,J(t,x;\alpha). $$
 The Hamilton-Jacobi-Bellman equation which allows to solve the problem above is (see \cite{4} )
\begin{equation}
\label{HJ}  \underset{\alpha\in U}{\rm sup}\Bigg( F(t,x;\alpha) +
\Big(\mathcal{L}^{\alpha}\Phi\Big)(t,x) \Bigg) =
0,
\end{equation}
where
\begin{equation*}
\Big(\mathcal{L}^{\alpha}\Phi\Big) =
\frac{\partial \Phi}{\partial t} + \alpha\frac{\partial
\Phi}{\partial x}
 + \frac{\delta^2}{2}\frac{\partial^2 \Phi}{\partial x^2},
\end{equation*}
with the terminal condition
\begin{equation}\nonumber
\label{(4)} \Phi(T,x) = K(x).
\end{equation}
 We denote the probability density of process  \eqref{X} 
 as $m(t,x)$.

 In what follows, we consider  a particular case
\begin{equation}\nonumber
\label{(1)}F(s,X_s;\alpha)= -\frac{\alpha^2}{2} + g(m(s,X_s),s,X_s).
\end{equation}
Thus, we assume that each agent receives a penalty for changing its position in the phase space and seeks to maximize its individual utility function, based on the fact that he knows only the distribution of the other agents. The choice of the function  $g$ depends on the type of problem.

Based on the specific form of the quality function  $F$, we obtain from
\eqref{HJ}
\begin{equation}\nonumber
\label{(5)}
\partial_t\Phi + \underset{\alpha}{\rm sup} \Big(\alpha\partial_x\Phi - \frac{1}{2}\alpha^2\Big) +
 g(m,t,x) + \frac{\delta^2}{2} \partial^2_{xx}\Phi = 0.
\end{equation}
Thus, if initially the density function $m(t,x)$ is known, we get the following  initial-terminal problem for  coupled Hamilton-Jacobi and Kolmogorov-Fokker-Planck equations
($x\in \mathbb R$, $t\in [0,T]$):
\begin{align}
\label{(6)}
&\partial_t\Phi + \frac{1}{2}(\partial_x\Phi)^2 +
\frac{\delta^2}{2}\partial^2_{xx}\Phi = -g,
\\\label{(7)}
&\partial_tm + \partial_x(m\,\partial_x\Phi ) =
\frac{\delta^2}{2}\partial^2_{xx} m,
\\\label{(8)}
&m(0,x) = m_0(x),\quad \Phi(x,T) = K(x).
\end{align}
It is known that under certain natural assumptions
(boundedness and Lipschitz continuity of  $g$
 and $K$)  problem \eqref{(6)}--\eqref{(8)} has only one classical solution  $(\Phi,m)$ \cite{5}.

The aim of this work is

1. 	To show that for some choice of $g$ and $K$ the solution of problem
\eqref{(6)}--\eqref{(8)} can be reduced to the solution of a system of Riccati equations.

2. 	To show that waiving the requirement of boundedness of $g$
 and $K$, we can construct an example  of nonexistence  of solution to
 \eqref{(6)}--\eqref{(8)}.

3. To show that the position of the maximum of density  $m(t,x)$ can be analytically found even though it is impossible to find solution of \eqref{(6)}--\eqref{(8)} for all time interval  $[0,T]$.

4.	To give an example of economic problem where  $g$ and $K$ have a form that allows to find a solution using the Riccati equations.

\section{Reduction to Riccati equations}
First of all, we note that from a mathematical point of view, it is interesting to analyze the function $ g $, which depends on $ m $, because in this case the equations \eqref {(6)} and \eqref {(7)} turn out to be linked. However, from the point of view of applications, the presence of $ m $ in the utility function of agents is not critical. Indeed, if  $g(m)$
increases then this means that it is beneficial for agents to stay closer to the maximum of  $m$. However, the presence of the term
$-\frac{\alpha^2}{2}$ already means that agents tend to resemble each other, this creates movement in the same direction.

Therefore, we consider a simpler case when $ g $ depends only on $ x $. We need it for applications below. We could, without prejudice to the method, add the term $\ln m$ to $g$  (which was done in [6]), however, for the sake of simplicity, we refuse this.

We also could consider a multidimensional analog of the problem ($x\in {\mathbb
R}^n$) via reducing the problem to the Riccati matrix equations. However,  the analysis of the obtained equations is a separate difficult problem and we strive to get the simplest model.

\subsection{Gaussian distribution}\label{Sec.2.1}

Let  $m_0(x) = M e^{-\frac{(x-x_0)^2}{\lambda}}$, where $\lambda>0$,
and $M>0$ is the normalization constant. Let also $g = a x^2 + bx +
c$, where $a, b, c$  are arbitrary constants. We  look for a solution to  system  \eqref{(6)}, \eqref{(7)}  in the form
\begin{equation}\nonumber
\label{f2} \Phi = A(t)x^2 + B(t)x + C(t), \quad
   m = e^{K_2(t)x^2 + K_1(t)x + K_0(t)},
\end{equation}
which imposes the terminal condition  $K(x)=A_T x^2+B_T x +C_T,$
with constants  $A_T, B_T, C_T$.

 Substitution into the system and equating the coefficients at the same powers of  $x$ gives the following set of equations:
\begin{align}
\label{f3} &A^{'} + 2A^2 = -a\\
\label{f4} &B^{'} + 2AB = -b,\\ &C^{'} + \delta^2 A^2+ \frac{B^2}{2}
= -c,\nonumber
\\
\label{f6}
&K_2^{'} + 4AK_2 - 2\delta^2K^2_2 = 0,\\
\label{f7} & K_1^{'} + 2AK_1 - 2\delta^2K_1K_2 + 2BK_2 = 0,
\\\nonumber
\label{f8} &K_0^{'} + 2A + BK_1 - \frac{\delta^2}{2}(K_1^2 + 2K_2) =
0.
\end{align}
The initial-terminal conditions are the following:
\begin{equation} \nonumber
K_2(0)=-\frac{1}{\lambda},\quad K_1(0)=-\frac{2 x_0}{\lambda},\quad
K_0(0)=-\frac{x_0^2}{\lambda}+\ln M,
\end{equation}
\begin{equation} \nonumber
A(T)=A_T,\quad B(T)=B_T,\quad C(T)=C_T.
\end{equation}
First, we note that equation \eqref{f3} which can be  elementary solved under the terminal condition $A(T)=A_T$ is separated from the system and determines its dynamics. Indeed, knowing  $A$, we find $K_2$ from
\eqref{f6}. Then we can find all the remaining functions.  It is easy to calculate that the solution is given by the following formulae:
\begin{itemize}
\item at $a<0$:
\begin{eqnarray}\label{A-}
  A(t) =-{k_-}\frac{(A_T-k_-){\rm e}^
{2\sqrt{-2a}(T-t)}+(A_T+k_-)}{(A_T-k_-) {\rm
e}^{2\sqrt{-2a}(T-t)}-(A_T+k_-)},\quad k_-=\sqrt{-\frac{a}{2}},
\end{eqnarray}
\item at
$a=0$:
\begin{eqnarray}\label{A0}
A(t) =\frac{A_T}{1-(T-t)A_T},
\end{eqnarray}
\item at $a>0$:
\begin{eqnarray}\label{A+}
 A(t) = \frac{1}{k_+}{\tan} \Big({\arctan}
\Big(k_+ A_T \Big) + \sqrt{2a}(T-t)\Big),\quad
k_+=\sqrt{\frac{a}{2}}.
\end{eqnarray}
\end{itemize}
Analyzing  \eqref{A-} and \eqref{A0}, it is easy to see that for  $a\le 0$
the solution  with  terminal condition $A_T<k_-$ exists for
$t\in [0,T]$ for all  $T>0$, while for $a>0$ the solution with any terminal condition $A_T$ can exist only for
$T<\frac{1}{\sqrt{2a}}(\frac{\pi}{2}-{\arctan}\,  (k_+A_T))$. When
$a\le 0$ and $A_T>k_-$ the solution also does not exist for all $T$.

In fact, these examples show that problem  \eqref{(6)} --
\eqref{(8)} can be unsolvable for all  $T$ for unbounded $g$ and $K$.

Further, we note that the position of maximum of
$m(t,x)$, ), which can be found as $Q = -\frac{K_1}{2K_2}$, can be described by  equation  $$Q^{'} = -\frac{K_1^{'}K_2 -
K_2^{'}K_1}{2K_2^{2}}.$$ Using  \eqref{f6}, \eqref{f7}, we get
\begin{equation}
\label{f9}Q^{'} =  2QA+B,
\end{equation}
and
\begin{equation*}
\label{f91} Q^{''} =  -2aQ-b.
\end{equation*}
 Thus, the position of  maximum either experiences fluctuations
(at $a>0$), or tends to reach a certain constant value outside  boundary layers near $t=0$ and  $t=T$. We know the initial condition  $Q(0)$, while to determine the second constant of integration it is necessary to use explicit expressions for $A$ and $B$.

Note that the function  $B$  can be found as a solution to the corresponding linear equations and written in elementary functions, but in a very cumbersome way, while the dynamics of  $Q$ is very simple.

Let us analyze the qualitative behavior of the solution of  \eqref{f3},
\eqref{f4}, \eqref{f9} with boundary conditions
\begin{equation}\label{bc}
A(T)=A_T,\quad B(T)=B_T,\quad Q(0)=Q_0.\nonumber
\end{equation}
First, we note that for  $a<0$ the solution of  system  \eqref{f3},
\eqref{f4}, has a stable equilibrium  $A_*=-k_-,$ $B_*=-\frac{b}{2A_*}$, to which the solution converges exponentially  under any terminal conditions ensuring
 the continuation of the solution on the entire semi-axis
   $(-\infty, T]$ .
Therefore, if  $T$ is large enough, then the solution of \eqref{f9}
outside the boundary layers near  $t=0$ and $t=T$  differs little from   $Q_*=-\frac{B_*}{2A_*}=-\frac{b}{2a}$. We will give pictures for specific data in the next section when we analyze an economical  problem.

Let us note that the function  $Q$ can be found for all  $T$, although it satisfies a linear first-order equation, the coefficients of which can become unbounded. For  $a=0$ , the form  $Q$
is the simplest:
\begin{equation}\label{Q0}
Q(t)=-\frac{1}{2} b t^2 +\frac{T(A_TT-1)b-B_T}{2A_TT-1} t
+\frac{2A_T(T-t)-1}{2A_TT-1} Q_0.\nonumber
\end{equation}
We can see that the formula is also defined for positive $A_T$, the only restriction is the condition $A_T\ne \frac{1}{2T}$. In this case, the position of the maximum has no limit as  $T\to \infty$
and with an increase in  $T$, the  maximum  of $m$ deviates unboundedly from its initial position.

If $a>0$, the position of maximum is described as follows:
\begin{equation}\label{Q+}
Q(t)=-\frac{b}{2a} + (Q_0+\frac{b}{2a})\cos (\sqrt{2a}t) + c_1 \sin
(\sqrt{2a}t),\nonumber
\end{equation}
where
$$c_1=\frac{1}{2\sqrt{a}\cos{\theta}}\Big[(b+2aQ_0)\sin(\theta)+\frac{b-aB_T}{\sqrt{A_T^2+k_+^2}}{\rm
sgn}(\cos\theta) \Big],$$ $$\theta={\arctan} \Big(\frac{ A_T}{k_+}
\Big) + \sqrt{2a}T.$$
 The position of the maximum of density varies periodically around
$-\frac{b}{2a}$.

For $a<0$ the solution is
\begin{equation}\label{Q-}
Q(t)=\frac{e^{2k_-t}}{(A_T-k_-)e^{4k_-T}-(A_T+k_-)}\left(Q_0 F_1 +
\frac{B_T}{2} F_2 +\frac{b}{8 k_-^2} F_3\right),\nonumber
\end{equation}
where $F_i=F_i(t,k_-,A_T$, $i=1,2,3$ are functions which do not have singularities, specifically,
\begin{align}&F_1=&((A_T+k_-)+(A_T-k_-)e^{4k_-(T-t)}),\quad F_2=(e^{
-4k_-t}-1)e^{2k_-T},\nonumber \\ &F_3=&(1-e^{-2k_-t})\left[
(1+e^{-2k_-t})((A_T-k_-)+(A_T+k_-)e^{2k_-T})\right.\nonumber\\&  &\left.-2(e^{2k_-(T-t)}
(A_T-k_-)+(A_T+k_-)) \right].\nonumber
\end{align}

\subsection{Distribution on the semi-axis  $[0,\infty)$.}\label{Sec.2.2}
Let  $m_0(x) = M x e^{-\frac{\kappa}{2} x^2}$, where $\kappa>0$, and
$M=\kappa$ is the normalization constant. For the solution of \eqref{(6)}
to be symmetric about zero  along the entire axis, we require  $$g = a
x^2 + c,$$ where $a,  c$ are arbitrary constants. We  find a solution to system  \eqref{(6)}, \eqref{(7)} in the form
\begin{equation}
\label{f2v} \Phi = A(t)x^2 + C(t), \quad
   m =x e^{-\frac{K_2(t)}{2} x^2 + K_0(t)},
\end{equation}
with the terminal condition $K(x(T))=A_T x^2+C_T,$ with constants $A_T,
C_T$.  In fact, we are looking for a solution to \eqref{(6)}, \eqref{(7)} with the boundary condition $ \partial_x\Phi|_{x=0}=m|_{x=0}=0$.

Substituting \eqref{f2v} in \eqref{(6)}, \eqref{(7)} gives the following system:
\begin{align}
\label{f3v} &A^{'}+ 2A^2 = -a,\\
&C^{'} + \delta^2 A^2= -c,\nonumber
\\
&K_2^{'}+ 4AK_2 + \delta^2K^2_2=0,\nonumber\\
&K_0^{'}+4A  + 3 \frac{\delta^2}{2}K_2=0.\nonumber
\end{align}
The initial-terminal conditions have the following form
\begin{equation}\label{bcv}
K_2(0)=\kappa,\quad K_0(0)=\ln M,\quad A(T)=A_T,\quad
C(T)=C_T.\nonumber
\end{equation}
The  maximum  of density $m$ is at the point
$Q(t)=\frac{1}{\sqrt{K_2(t)}}$. It is easy to calculate that the position of the maximum satisfies the following Bernoulli equation:
\begin{equation}\label{Qv}
Q'=2AQ+\frac{\delta^2}{2Q},
\end{equation}
where $A$ can be found according \eqref{A-}, \eqref{A0} and
\eqref{A+}. The system  \eqref{f3v}, \eqref{Qv} has the first integral
$$(a+2A^2(t))Q^2+\delta^2A=\rm const,$$
which allows  to express  $Q$ from $A$, specifically:
\begin{equation}\label{Qvs}
Q(t)=\left( \frac{(a+2A^2(0))Q^2(0)-\delta^2(A(t)-A(0))}{a+2A^2(t)}\right)^\frac12, \,
\end{equation}
where $Q(0)=\kappa^{-1/2}$. From \eqref{Qvs} we can find that, as in  Sec.\ref{Sec.2.1}, the function $Q(t)$ does not have singularities inside the interval  $[0,T]$, even if $A(t)$
goes to infinity. In these points $Q$ tends to zero. If $a<0$,  then for large $T$ the maximum of density  $m(t,\xi)$ is close to the equilibrium
$Q_*=\frac{\delta}{4\sqrt{k_-}}$.

\section{An example of application of the mean field game theory: forming the opinion of investors about the asset}\label{Sec.3}

We give an example showing that the solutions found in the previous section have natural applications in the field of financial mathematics. To do this, consider a market in which a large number of investors operate, managing their own portfolio of securities, consisting of a risky asset and a deposit, by solving the Merton problem \cite{7}.

\subsection{Individual strategy}\label{Sec.3.0}
Let us retell the statement of the problem, which each investor solves individually, following \cite{6}, Example 11.2.5.

The price of the risky asset
$S_1$ is described by the stochastic differential equation
\begin{equation}
\label{(77)} dS_1 = \mu S_1dt + \sigma S_1dW_t,
\end{equation}
where $W_t$  is the standard Wiener process $\mu=\rm const$,
$\sigma={\rm const}
> 0$. In the economic context, these values are commonly called the drift parameter and the volatility parameter, respectively.
The price of the risk-free asset $S_2$ is determined only by a constant interest rate $r$:
\begin{equation} \label{(8.1)} dS_2 = r S_2dt.\nonumber
\end{equation}
Suppose that the investor operates the portfolio  $V$, which consists of risky and risk-free assets, with $ h_1 $ and $ h_2 $ being the shares of capital invested in risk and risk-free assets, respectively, $h_1 + h_2 = 1$.
  Then
\begin{equation}
\label{(9.1)} \frac{dV}{V} = h_1\frac{dS_1}{S_1} +
h_2\frac{dS_2}{S_2}.\nonumber
\end{equation}
Let us denote $ h_1=h$, $h_2=1-h$. The change in the value of the portfolio  has the following form::
\begin{equation}
\label{(V)} dV= (r + (\mu - r)h)V\,dt + \sigma h V\,  dW_t.
\end{equation}
Assume that starting with capital $ V_t = v> 0 $ at the time $ t $, the investor wants to maximize the expected return on capital at some subsequent point of time $ T> t $.
 If we set the utility function $ N (V) $, which is usually assumed to be increasing and convex upward, then the problem reduces to finding the function $ \Phi (t, v) $ and the Markov control $ h_* = h _* (t, V) $, such  that
\begin{equation}
\label{(11p)}\Phi(t,v)=\underset{h}\sup \{J^h(t,v)\}=
J^{h_*}(t,v),\nonumber
\end{equation}
where $h$ is the Markov control, $J^h(t,v)={\mathbb
E}^{t,v}[N(V_T^h)]$. In order to solve this problem we should define a differential operator
\begin{equation}
\label{(13)} \mathcal{L}^h f = \frac{\partial f}{\partial t} + (\mu
h+r (1-h))v \,\frac{\partial f}{\partial v} + \frac{1}{2}\sigma^2
h^2 v^2 \,\frac{\partial^2 f}{\partial v^2}\nonumber
\end{equation}
and solve the Hamilton-Jacobi-Bellman equation:
\begin{align}
\label{HJB}\sup_h \{({\mathcal L}^h \Phi)(t,v)\}=0,\quad t\in(0,T),
\,v>0,\\
\Phi(T,v)=N(v),\quad \Phi(t,0)=N(0), \quad t<T.\nonumber
\end{align}
If $\partial_v\Phi>0$ and $\partial^2_{vv}\Phi<0$,  then the solution is
 $h(t,v)=\frac{(\mu-r)\partial_{v}\Phi}{v\sigma^2\partial^2
_{vv}\Phi}$. Substituting this expression in \eqref{HJB} gives
the following boundary value problem for  $\Phi$:
\begin{align}
\label{HJB1}
\partial_t\Phi+rv\partial_v\Phi-\frac{(\mu-r)\partial_v\Phi}{2\sigma^2\partial^2_{vv}\Phi}=0, \quad
t\in(0,T),
\,\,v>0,\\
\Phi(t,v)=N(v), \quad t=T \,\,\mbox{or}\,\, v=0.\nonumber
\end{align}

As a utility function we take $ N (v) = \frac{v^q}{q}$, $ q <1 $, $ q \ne 0 $, or $N(v)=\ln v$.
The latter function formally corresponds to the limit $ q \to 0 $. All these functions belong to the class HARA (hyperbolic absolute risk aversion) \cite{8}.
 In addition, the case $ q <0 $ corresponds to the strategy of the investor who prefers the least risky investments, $ q = 0 $ corresponds to risk-neutral strategies, $ q> 0 $ corresponds to the risk-prone investor \cite{9}, Sec.2.

The solution \eqref{HJB1} can be found in the form $\Phi(t,v) = \phi(t)N(v)$,
the corresponding optimal strategy at all $q<1$ is
\begin{align}\label{OS}
 h_* = \frac{(\mu - r)}{\sigma^2(1-q)}.
\end{align}

Let us calculate the capital  growth rate  $\displaystyle{\mathbb E}\frac{\ln
V}{t}$, that investor will receive guided by the optimal strategy \eqref{OS}. From \eqref{(V)} and Ito's formula:
\begin{equation}
\label{(22)} d \ln V = \Big[(\mu h + r(1-h)) - \frac{1}{2}\sigma^2
h^2 \Big]dt + \sigma h \,dW_t,\nonumber
\end{equation}
for all $q<1$
\begin{equation}
\label{(24)}{\mathbb E}\frac{\ln V}{t}= r +
\frac{(1-2q)(\mu-r)^2}{2\sigma^2(q-1)^2}.\nonumber
\end{equation}
This, in particular, implies that strategies for investing in risky assets with  $q>\frac12$  (very risky investors) lead to a decrease in portfolio returns.

\subsection{Collective strategy}\label{Sec.3.1}
Now let us describe the collective behavior of investors.
We assume that they all manage the portfolio based on their own ideas about the parameters (drift and volatility parameters) of the risk asset. In other words, each fixed investor carries out
control based on the equation \eqref {(77)} with its own choice of $ \mu $ and $ \sigma $.
The УtrueФ values of these parameters
(we denote them by  $\bar\mu$ and $\bar\sigma$) are unknown to investors.
 These УtrueФ values may differ from those accepted in the market, and they manifest themselves only in that the investor receives a penalty for their incorrect choice.

It is believed that the opinions of investors about the correct value of the drift parameter  $\mu$ are distributed normally along the entire axis, the maximum is initially at $\mu_0$.

 The volatility $\sigma$
obeys  some positive distribution and its density initially has  a maximum at some point $\sigma_0>0$.
 At the same time, investors receive a penalty both for deviating from the "true" values of the drift and volatility parameters, and for deviating from the majority opinion.
A significant simplification that stems from the desire to obtain an analytical solution to the problem is the assumption that all investors treat risk the same way, guided by the same utility function $ N(v)$.

 During the control process, the initial distributions of drift and volatility parameters change with the desire to maximize capital growth rate.

We will be interested in how the position of the distribution maximum  $m(t,x)$ changes in response to the control method, that is, how the market is forming an opinion about the parameters of a risky asset.

We obtain a typical optimization problem of the theory of mean-field games
\eqref{(1)}, when the random variables $\mu$ or
$\phi(\sigma)$ play the role of $X$, subordinate to \eqref{X}, $K(X(T))= \frac{\ln
V(T)}{T}$, and
\begin{equation}
g(t,X, m)=\beta\frac{\ln V(X)}{t} - \gamma\eta(X,\bar{x}) +  \lambda
\ln\frac{m(t,X)}{\mu_*(t)},\nonumber
\end{equation}
$\beta\ge 0$, $\gamma\ge 0$, $\lambda\ge 0$,   $ \bar{x}=\rm const$,
$\phi(\sigma)$  is some smooth function of  volatility, $\eta$
is a function simulating  penalty for an investor who incorrectly guesses the parameters of the risky asset. It is chosen for reasons of convenience, for example, $\eta(X,\bar{x})=(X - \bar{x})^2$. The presence of the function $ \ln\frac{m}{\mu_*}\le 0$,
 where $\mu_*(t)>0$  is the value of the maximum of density $m(t,x)$ at a fixed time  $t$
models a penalty for deviating the opinion of a investor about a risky asset from the majority opinion.
 As we have already noted, the first term under the integral in \eqref{(1)} serves the same purpose. Moreover, it can be analytically shown that the equation describing the position of the distribution maximum does not change if we assume $\lambda>0$, the presence of this term will only lead to a more pronounced maximum, which is also confirmed by numerical analysis. Therefore, below we assume $\lambda=0$.
The coefficients  $\beta$ and $\gamma$  can be considered equal to zero or not, depending on what problems we are studying.

As application of the results of the previous sections, we consider two separate cases.  In the first of them, we suppose that the volatility of a risky asset is known and an opinion is formed regarding its drift. In the second, on the contrary, the drift parameter is considered known, and an opinion is formed about volatility.

\subsection{The opinion about the parameter of drift $\mu$}\label{Sec.3.2}
So, we assume that the volatility  $\sigma$ is fixed,
$\eta=(\mu - \bar{\mu})^2$. The system of equations
\eqref{(6)}--\eqref{(8)} has the following form:
\begin{align}
\label{(33)} &\partial_t\Phi + \frac{1}{2}(\partial_\mu\Phi)^2 +
\frac{\delta^2}{2}\partial^2_{\mu \mu}\Phi = -\beta (r +R(\mu -
r)^2) + \gamma(\mu - \bar{\mu})^2,
\\
&\partial_t m + \partial_{\mu}(m\partial_\mu\Phi) -
\frac{\delta^2}{2}\partial^2_{\mu \mu} m = 0,\nonumber
\\
&\Phi(\mu, T) = r + R(\mu - r)^2, \quad m(\mu, 0) =
\frac{1}{\sqrt{\pi \lambda}}e^{-\frac{(\mu-\mu_0)^2}{\lambda}},
 \label{(33.1)}
\end{align}
where $R= \frac{(1-2q)}{2\sigma^2(q-1)^2}$.

Thus, in the notation of Sec.\ref{Sec.2.1}
\begin{align}
&a=\beta R-\gamma,\quad b=2(\gamma\bar \mu-\beta Rr),\quad c=\beta r+\beta Rr^2-\gamma\bar \mu^2,\nonumber\\
& K_2(0)=-\frac{1}{\lambda},\quad K_1(0)=-\frac{2
\mu_0}{\lambda},\quad K_0(0)=-\frac{\mu_0^2}{\lambda}+\ln
\frac{1}{\sqrt {\pi\lambda}},\nonumber\\&A_T=R,\quad B_T=-2Rr,\quad
C_T=r+Rr^2.\nonumber
\end{align}

According to the results of Sec.\ref{Sec.2.1}, the solution of
\eqref{(33)}--\eqref{(33.1)} exists for all $T>0$, if and only if $\beta R-\gamma<0$ and $R<\sqrt{\frac{\gamma-\beta R}{2}}$.
However, the position of the maximum of  $m(t,\mu)$ can always be determined. Namely, if $\beta R-\gamma<0$ ($a<0$), then for large $T$ the position of maximum  is close to  $$Q_*=\frac{r\beta
R-\gamma\bar \mu}{\beta R-\gamma}.$$ In other words, in this case, investors form an opinion about the asset. If $\beta
R-\gamma>0$ ($a>0$), then the maximum of  $m$ oscillates periodically, sometimes  deviating significantly from its average value. The frequency of these oscillations increases with  $a$. In this case, we say that investors cannot agree on the parameters of asset.

Let us analyze this result from the point of view of the agent's behavior when investing for a long period of time.

\begin{itemize}

\item If $R<0$, that is, investors adhere to a rather risky strategy ($q>\frac12$), then the market forms the opinion about the correct value of $\mu$, even if  $\gamma=0$. In the threshold case $q\to 1$ most investors believe that the correct return on a risky asset is close to the risk-free rate of the asset and avoid investing in a risky asset. If $q\to\frac12$, then, otherwise, the opinion tends to the "true" value of $\mu$.

\item  If $R>0$, that is, the investors are rather cautious
($q<\frac12$), then
\begin{itemize}
\item  without the presence of a penalty for the wrong choice of $\mu$
an opinion on the return of the asset does not form.

\item If $\beta=0$ or
$\gamma$ is sufficiently large, then  $Q_*$ is close to the "true" value
$\mu=\bar \mu$.

\item The similar effect has an unlimited increase of
$\sigma$ or $|q|$, since $\lim\limits_{\sigma \to +\infty} R
=\lim\limits_{q\to -\infty} R=0$.

In other words, very cautious investors quickly form the correct opinion about the asset, and this
the faster, the greater the volatility.
\end{itemize}

 \item	Since  $\lim\limits_{\sigma \to 0} R =\infty$, the smaller the randomness component in the value of the risky asset, the more difficult it is for investors to come to a common opinion about it.
     Indeed, to ensure the condition $\beta R-\gamma<0$  for a fixed $q$, so we should  choose a large value of $\gamma$.  This situation seems paradoxical, but can be explained as follows: if the investor is dealing with a low-risk asset, then he is inclined to adhere to a more risky strategy, that is, choose a larger $ q $. It is easy to see that if
$q=\frac12-o(\sigma^2)$, then $R$ tends to zero, which helps to determine the correct value of $\mu$.

\end{itemize}

\begin{figure}[h]\label{fig1}
\centerline{\includegraphics[width=0.5\columnwidth]{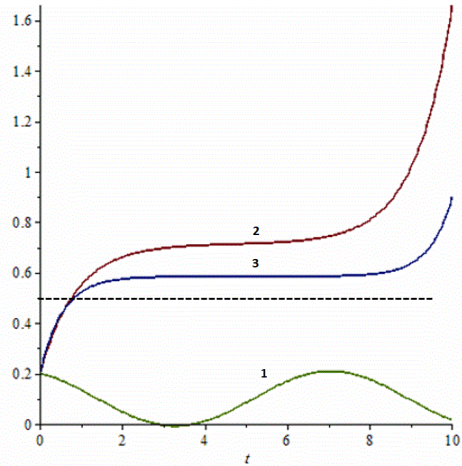}}
\caption{The position of  maximum of $\mu$ for $q=-10$,
$\sigma=0.5$, $\beta=1$, $\bar \mu=0.5$, $r=0.1$ at values
$\gamma=0, \, 1,\, 2$ (marked by numbers). $Q_0=0.2$. With increase of
$\gamma$ the  maximum  tends to the correct value, with a tendency to overestimate.}
\end{figure}

The typical behavior of the maximum of density $m(t,\mu)$ is represented on Fig.\ref{fig1}.

\subsection{Opinion on volatility $\sigma$}\label{Sec.3.3}

Now suppose that the parameter $\mu$ of the asset is known, and the agents make their assumptions regarding volatility. Unlike the previous problem, we cannot assume that $\sigma$ is distributed on the entire axis and we should consider some positive distribution of it. In order to be able to obtain an analytical solution, we consider instead of $\sigma$ the quantity  $\xi =
\frac{1}{\sigma}$, which at the initial moment of time has the density $m_0(\xi) = \kappa \xi e^{-\frac{\kappa}{2} \xi^2}$
with a maximum at  $\xi_0=\frac{1}{\sigma_0}=\frac{1}{\sqrt{\kappa}}$ and  use the results of Sec.\ref{Sec.2.2}. Due to the fact that the function  $g $ must be even in $\xi$, we cannot choose $\eta$ such as in Sec.\ref{Sec.3.2}. Instead, we assume that $\eta
(\xi,\bar\xi)=\xi^2$, that is, the investor receives a penalty for choosing $\xi$ too large (or $\sigma$ too small),
formally  $\bar\xi=0$. We emphasize that we change the problem for the sake of the ability to get an exact solution.

The problem \eqref{(6)}--\eqref{(8)} has the following form
\begin{align}
& \partial_t\Phi + \frac{1}{2}(\partial_{\xi}\Phi)^2 +
\frac{\delta^2}{2}\partial^2_{\xi \xi}\Phi = -\beta(r + P\xi^2) +
\gamma \xi^2,\nonumber
\\
& \partial_t m + \partial_{\xi}(m \partial_{\xi} u) -
\frac{\delta^2}{2}\partial^2_{\xi \xi}m =
0,\nonumber\\
&\Phi(\xi, T) = r + P \xi^2,\quad m(\xi, 0) = \frac{1}{\xi_0^2 }\xi
e^{-\frac{1}{2\xi_0^2} \xi^2},\nonumber\label{(33.1)}
\end{align}
where $P= \frac{(1-2q)(\mu - r)^2}{2(q-1)^2}$. In notation of  Sec.\ref{Sec.2.2}
\begin{align*}
a=\beta P-\gamma,\quad c=\beta r, \end{align*}
\begin{align*}
K_2(0)=\frac{1}{\xi_0^2},\quad K_0(0)=\ln \frac{1}{\xi_0^2},\quad
A_T=P,\quad C_T=r.
\end{align*}

Let us make conclusions about the behavior of the maximum of $m(t,\xi)$,
which which follow from the results of Sec.\ref{Sec.2.2}.

\begin{itemize}

\item If we assume that  $\gamma=0$, $\beta>0$ then for large $T$
the maximum $Q_*$ is close to a constant if and only if $P<0$, that is
$q>\frac{1}{2}$.  This constant is equal to  $\frac{\delta}{\sqrt{-8 P}}$.
We associate this behavior with forming an opinion on the volatility of a risky asset. The solution reaches a constant, the faster, the greater the difference $\mu-r $. If $\delta$ decreases, the value of  $Q_*$ tends to zero. This means that the  maximum of $m(t,\sigma)$ tends to infinity, that is, the risky asset is perceived as more uncertain.

\item If $\gamma\ge 0$ and $P>0$, then the opinion on volatility is formed under the condition  $\beta P-\gamma<0$. For large $T$
the maximum $Q_*$ is close to a constant equal to $\frac{\delta}{\sqrt{8(\gamma-\beta P)}}$.

\item If $\beta P-\gamma>0$, then the position of the  maximum of
$m(t,\sigma)$ oscillates with a frequency that grows with  $\beta P-\gamma$.

At a qualitative level, the situation is identical to that described in Sec.
\ref{Sec.3.2}.

\end{itemize}

\section{Conclusion}
We consider a simple application of the theory of mean field games to study  the behavior of market agents managing a portfolio of securities which consist of risky and risk-free assets, based on a utility function that is common to all.  We assume that the information on the market is incomplete, that is, agents are forced to independently decide on the parameters of the risky asset. When setting problems, we limit ourselves to the possibility of obtaining its analytical solution.

  We deal with two separate situations. In the first one, agents know the exact value of volatility parameter, but decide on the correct value of the drift parameter. When agents manage their portfolio, they get a penalty for a false assumption of the "true" value of the drift parameter. In the second case, on the contrary, the value of the drift parameter is known, however, agents receive a penalty for considering volatility too small. We study the question of whether, under the conditions described, the market  formed an opinion about parameters of the asset, and if so, how far is it from the correct one.

 The model that can be extended in different directions. In particular, it is natural to assume that agents have to choose three parameters at the same time: the drift and volatility of the asset, as well as risk attitude. Such problem can also be solved within the framework of the mean-field games theory, however, it is three-dimensional in space and does not allow an analytical solution. However, it can be investigated numerically.

\end{document}